\documentclass{article}

\usepackage[final]{neurips_2020}

\usepackage[utf8]{inputenc} 
\usepackage[T1]{fontenc}    
\usepackage{hyperref}       
\usepackage{url}            
\usepackage{booktabs}       
\usepackage{amsfonts, amsmath}       
\usepackage{nicefrac}       
\usepackage{microtype}      

\usepackage[utf8]{inputenc}
\usepackage{amsmath,amssymb,hyperref,array,xcolor,multicol,verbatim, dsfont}
\usepackage[normalem]{ulem}
\usepackage[pdftex]{graphicx}
\usepackage{fullpage}
\usepackage{enumitem}
\usepackage{algorithm}
\usepackage{algorithmic}
\usepackage{natbib}
\hypersetup{
	pdfauthor={Anonymous},     
	pdfsubject={},   
	pdfcreator={Anonymous},   
	pdfproducer={},  
}
\setenumerate[1]{label=(\arabic*)}
\newtheorem{remark}{Remark}
\newtheorem{theorem}{Theorem}

\usepackage{chngpage}
\newtheorem{example}{Example}
\usepackage{lipsum}
\newcommand{\eins}{\text{\ensuremath{1\hspace*{-0.9ex}1}}}

\title{MinMax Methods for Optimal Transport and Beyond: Regularization, Approximation and Numerics}

\author{%
	Luca De Gennaro Aquino\thanks{Equal contribution}\\
	Department of Accounting, Law \& Finance\\
	Grenoble Ecole de Management\\
	38000 Grenoble, France\\
	\texttt{luca.degennaroaquino@grenoble-em.com} 
	\And
	Stephan Eckstein$^{*}$\\
	Department of Mathematics and Statistics\\
	University of Konstanz\\
	78464 Konstanz, Germany\\
	\texttt{stephan.eckstein@uni-konstanz.de}
}

\begin{document}
	\maketitle
	\begin{abstract}
		We study MinMax solution methods for a general class of optimization problems related to (and including) optimal transport. Theoretically, the focus is on fitting a large class of problems into a single MinMax framework and generalizing regularization techniques known from classical optimal transport. We show that regularization techniques justify the utilization of neural networks to solve such problems by proving approximation theorems and illustrating fundamental issues if no regularization is used. We further study the relation to the literature on generative adversarial nets, and analyze which algorithmic techniques used therein are particularly suitable to the class of problems studied in this paper. Several numerical experiments showcase the generality of the setting and highlight which theoretical insights are most beneficial in practice.
	\end{abstract}
	\section{Introduction}
	Optimal transport (OT) has received remarkable interest in recent years (\cite{cuturi2013sinkhorn, peyre2019computational}). In many areas, classical OT and related problems with either additional constraints (\cite{tan2013optimal, korman2015optimal, beiglbock2016problem, nutz2020directional}) or slight variations (\cite{buttazzo2012optimal,pass2015multi,Chizat2018UnbalancedOT, liero2018optimal, seguy2017large}) have found significant applications. In this paper, we propose a MinMax setting which can be used to solve OT problems and many of its extensions and variations numerically. In particular, the proposed methods aim at solving even non-discrete problems (i.e., where continuous distributions occur), which is often relevant particularly for applications within finance and physics.
	
	The basic premise of solving the described class of problems with MinMax techniques and neural networks has been applied in less general settings various times (\cite{yang2018scalable,xie2019scalable,henry2019martingale}). 
	The idea is that one network generates candidate solutions to the optimization problem, and a different network punishes the generating network if the proposed candidates do not satisfy the constraints.
	Both networks play a zero-sum game, which describes the MinMax problem (see Section \ref{sec:theory} for more details).
	Compared to the widely known numerical methods based on entropic regularization (\cite{cuturi2013sinkhorn, Solomon2015ConvolutionalWD, genevay2016stochastic}), where measures are represented via their densities, in MinMax settings the candidate solution is expressed as a push-forward measure of a latent measure under the function represented by the generator network. 
	
	Within the class of MinMax problems studied in the literature, some instances are easier to solve numerically than others. Particularly for classical OT problems, regularization techniques that slightly change the given optimization problem, but lead to better theoretical properties or numerical feasibility, have been proposed (\cite{yang2018scalable,xie2019scalable}).  Section \ref{sec:reformulations} showcases how such regularization techniques can be applied in the general framework at hand. 
	
	While the theoretical approximation by neural networks of the general MinMax problem  can fail even in very simple situations if no regularization is used (see Remark \ref{rem:pmapprox}), regularization techniques can yield substantial improvements in this regard (see Theorem \ref{thm:approxreg}). We emphasize that it is of fundamental importance to couple the implemented problem utilizing neural networks to theory. In classical (non MinMax) optimization problems, the ubiquity of universal approximation theorems may give the false impression that any approximation by neural networks within optimization problems is justified. It is a key insight of this paper that in MinMax settings this becomes more difficult, but can still be achieved with the right modeling tools (e.g., regularization of the theoretical problems). In relation to generative adversarial networks (GANs) (\cite{goodfellow2014generative}), the results established in \mbox{Section \ref{sec:reformulations}} can be regarded as new insights as to why regularization is helpful during training (see, e.g., \cite{roth2017stabilizing}), and we show that this improved stability can also be observed in the numerical experiments in Section \ref{sec:numerics} and Appendix \ref{app:w2}.
	
	With or without regularization techniques, MinMax problems utilizing neural networks remain notoriously hard to solve numerically. Within the literature on GANs, this has sparked several algorithmic and numerical ideas that go beyond simple gradient descent-ascent (GDA) methods. Section \ref{sec:algorithmic} discusses what we consider to be the most relevant methods that should be adapted to the class of problems studied in this paper. 
	
Finally, Section \ref{sec:numerics} reports numerical experiments. The experiments showcase how the insights obtained by the theoretical regularization techniques, and the relation to the GAN literature, can be utilized for numerical purposes. To illustrate the generality of the setting, we go beyond classical OT problems by taking examples from optimal transport with additional distribution constraints and martingale optimal transport. The Appendix includes technical proofs, details about the numerical experiments, a list of problems from the literature and further numerical results.


	\section{Theoretical setting}
	\label{sec:theory}
	Let $d\in \mathbb{N}$, $\mathcal{P}(\mathbb{R}^d)$ be the Borel probability measures on $\mathbb{R}^d$, and $C(\mathbb{R}^d)$ (resp. $C_b(\mathbb{R}^d)$) be the space of continuous (and bounded) functions mapping from $\mathbb{R}^d$ to $\mathbb{R}$. 
	Let $\mu \in \mathcal{P}(\mathbb{R}^d)$, $f\in C(\mathbb{R}^d)$ and $\mathcal{H} \subset C(\mathbb{R}^d)$.
	The optimization problem studied in this paper is of the form 
	\begin{equation}
	\label{eq:primaldef}
	(P) = \sup_{\nu \in \mathcal{Q}} \int f \,d\nu, \text{ where }
	\mathcal{Q} := \left\{ \nu \in \mathcal{P}(\mathbb{R}^d) : \int h \,d\nu = \int h \,d\mu \text{ for all } h \in \mathcal{H}\right\}.
	\end{equation}
	The class of problems of the form $(P)$ can be seen as 
	the class of 
	linearly constrained problems over sets of probability measures. Most 
	relevant is certainly the subclass of linearly constrained 
	optimal transport problems (\cite{zaev2015monge}). For an incomplete but 
	illustrative list of further examples from the literature we refer to 
	Appendix \ref{subsec:listofproblems}.
	The most popular representative within this class of problems is the following:
	\begin{example}[Optimal transport]
		\label{ex1}
		Let $d = 2$ and $\mathcal{H} = \{h \in C_b(\mathbb{R}^2) : \exists \, h_1, h_2 \in C_b(\mathbb{R}): \forall (x_1, x_2) \in \mathbb{R}^2 : h(x_1, x_2) = h_1(x_1) + h_2(x_2)\}$ and let $\mu \in \mathcal{P}(\mathbb{R}^2)$. Then it holds
		\[\mathcal{Q} = \left\{ \nu \in \mathcal{P}(\mathbb{R}^2) : \nu_1 = \mu_1 \text{ and } \nu_2 = \mu_2\right\} =: \Pi(\mu_1, \mu_2),\]
		where $\mu_j$ and $\nu_j$, for $j=1, 2$, denote the $j$-th marginal distribution of $\mu$ and $\nu$, respectively.
	\end{example}
	The example shows that the precise choice of $\mu$ is often irrelevant, and 
	only certain characteristics of $\mu$ (like its marginal distributions in 
	the case of optimal transport) are relevant.
	
	To make the set $\mathcal{H}$ more explicit and allow for an approximation by neural networks, we restrict the form of $\mathcal{H}$ to
	\begin{equation}
	\mathcal{H} = \left\{ \sum_{j=1}^J e_j \cdot (h_j \circ \pi_j) : h_j \in C_b(\mathbb{R}^{d_j})\right\},
	\label{FH}
	\end{equation}
	where $J\in \mathbb{N}$ and $e_j : \mathbb{R}^d \rightarrow \mathbb{R}$ and $\pi_j: \mathbb{R}^d \rightarrow \mathbb{R}^{d_j},$ for all $j=1,\dots, J,$ are fixed. 
	This form of $\mathcal{H}$ is not a strong restriction, as it includes all 
	relevant cases that the authors are aware of. Hereby, the 
	functions $\pi_j$ 
	can be seen as transformations of the input variable. For instance, the projection
	$\pi_j(x) = x_j$ onto the $j$-th variable is used 
	in optimal transport. The functions $e_j$ are used to scale the transformed input. While for OT $e_j \equiv 1$, for instance in the MOT problem 
	(\cite{beiglbock2013model}) one sets
	$e_j(x_1, x_2) = x_2 - x_1$ to enforce the martingale constraint.
	
	This form of $\mathcal{H}$ now allows for a neural network approximation. We fix a continuous activation function and a number of hidden layers. For $d_1, d_2, m \in \mathbb{N}$, let $N_{d_1, d_2}^m$ be the set of all feed-forward neural network functions mapping $\mathbb{R}^{d_1}$ to $\mathbb{R}^{d_2}$ with hidden dimension $m$ (by hidden dimension we mean the number of neurons per layer).
	We then define the neural network approximation $\mathcal{H}^m$ of $\mathcal{H}$ by
	\begin{equation}
	\label{FHm}
	\mathcal{H}^m := \left\{ \sum_{j=1}^J e_j \cdot (h_j \circ \pi_j) : h_j \in N^m_{d_j, 1}\right\}.
	\end{equation}

	\subsection{Reformulation as MinMax problem over neural network functions}
	This subsection shows how to reformulate $(P)$ as an \textsl{unconstrained} optimization problem over neural network functions, which leads to a MinMax problem.
	Let $K \in \mathbb{N}$ and $\theta \in \mathcal{P}(\mathbb{R}^{K})$. For a function $T$, denote by $\theta_{T} := T_* \theta$ the push-forward measure of $\theta$ under $T$. The following approach builds on representing arbitrary probability measures $\nu \in \mathcal{P}(\mathbb{R}^d)$ as the push-forward of $\theta$ under some measurable map $T$. To make this work, $\theta$ has to allow for sufficiently rich variability of its push-forward measures. This means, e.g., $\theta$ may not simply be a discrete distribution. More precisely, we can require that $\theta$ can be reduced to the uniform distribution on the unit interval $(0, 1)$, which suffices so that any measure $\nu \in \mathcal{P}(\mathbb{R}^d)$ can be written as $\nu = T_*\theta$ for some measurable map $T$.\footnote{Formally, the argument works as follows: Denote by $\mathcal{U}$ the uniform distribution on $(0,1)$. Say there exists a measurable map $S : \mathbb{R}^K \rightarrow (0,1)$ such that $\mathcal{U} = S_*\theta$. Let $\nu \in \mathcal{P}(\mathbb{R}^d)$ be arbitrary. We know that there exists a bimeasurable bijection $B: (0, 1) \rightarrow \mathbb{R}^d$. Let $B_{\textrm{inv}}$ denote the inverse of $B$. Set $\tilde{\nu} = (B_{\textrm{inv}})_* \nu$ and let $Q_{\tilde{\nu}}$ be the quantile function of $\tilde{\nu}$. Then $\nu = (B \circ Q_{\tilde{\nu}})_* \mathcal{U} = (B \circ Q_{\tilde{\nu}} \circ S)_* \theta$.} The measures $\theta$ used in practice are usually high-dimensional Gaussians or uniform distributions on unit cubes, which all satisfy this requirement. We reformulate $(P)$ as follows:
	\begin{equation}
		\begin{split}
			(P) & = \sup_{ \nu \in \mathcal{P}(\mathbb{R}^{d})} \inf_{h \in \mathcal{H}} \int f \, d\nu + \int h \, d\nu - \int h \, d\mu \\
			& = \sup_{T : \mathbb{R}^{K} \to \mathbb{R}^{d}} \inf_{h \in \mathcal{H}} \int f\, d\theta_T + \int h \, d\theta_T - \int h\, d\mu \\
			& = \sup_{T : \mathbb{R}^{K} \to \mathbb{R}^{d}} \inf_{h \in \mathcal{H}} \int f \circ T \, d\theta + \int h \circ T  \, d\theta - \int h\, d\mu \\	
			& \approx \sup_{T \in N^m_{K,d}} \inf_{h \in \mathcal{H}^m} \int f \circ T \,  d\theta + \int h \circ T\, d\theta - \int h \, d\mu =: (P^m) \label{MM} \\
		\end{split}
	\end{equation}
	
	In Eq.~\eqref{MM}, the approximation $\approx$ by neural networks is more subtle than a simple application of a universal approximation theorem. At many points in the literature, this approximation, and particularly the subtle difficulties that occur due to the sup-inf structure, have been overlooked. A related MinMax problem occurs in the Projection Robust Wasserstein distance (\cite{paty2019subspace,lin2020projection}) and similar approximations to $\approx$ above have been studied by \cite{champion2004gamma,degiovanni2014limit}.
	\begin{remark}
	\label{rem:pmapprox}
	\begin{itemize}
		\item[(i)] In general, the approximation by neural networks in Eq. \eqref{MM} fails, i.e., it might \textbf{not} hold $(P^m) \rightarrow (P)$ for $m\rightarrow \infty$. 
		
		A simple counterexample is the following optimal transport example (see Example \ref{ex1}): Let $\mu_1 \in \mathcal{P}(\mathbb{R})$ be given by its Lebesgue density $\frac{d\mu_1}{d\lambda}(x) := {\rm 1}_{(0, 1)}(x) \kappa x \cdot |\sin(1/x)|\,$ for a suitable $\kappa > 0$, and set $\mu_2 = \mu_1$. Let $\theta$ be the uniform distribution on $(0,1)$ and let the activation function of all networks be the ReLU function. Then $(\theta_T)_1 \neq \mu_1 \,$ for all $T \in N_{1, 1}^m \,$ and hence $\,\inf_{h \in \mathcal{H}^m} \int h \,d\theta_T - \int h \,d\mu = -\infty$, since $(x \mapsto a\cdot (x-b)^+) \in \mathcal{H}^m$ for all $a, b \in \mathbb{R}$.\footnote{The reason the infimum evaluates to $-\infty$ is that for two measures $\nu \neq \mu$ one can find $b \in \mathbb{R}$ such that $\int (x-b)^+\,\mu(dx) \neq \int (x-b)^+ \,\nu(dx)$ (see \cite[Footnote 2]{beiglbock2013model}).} Thus, if $f\equiv 0$, then $(P) = 0$, while $(P^m) = - \infty$ for all $m$.
		
		We see that it does not matter how closely $(\theta_T)_1$ approximates the marginal distribution $\mu_1$, as any deviation can be exploited arbitrarily by the inner infimum problem. Both regularization techniques for problem $(P)$ introduced in Section \ref{sec:reformulations} will resolve this issue.
		
		Further, the theme of the above counterexample is quite general. Whenever a distribution is precisely specified (like a marginal distribution in optimal transport) and it cannot be represented exactly by the network $T$, then the inner infimum will evaluate to minus infinity.
		
		\item[(ii)] The problem in (i) is not that neural networks lack approximation capabilities. Indeed, define $\Phi(T, h) := \int f \circ T \, d\theta + \int h \circ T  \, d\theta - \int h\, d\mu$. In the setting of (i), and more generally, applying standard universal approximation theorems can show that for any $T:\mathbb{R}^K \rightarrow \mathbb{R}^d$ Borel and $h \in \mathcal{H}$, there exist $T^m \in N_{K, d}^m$, $h^m \in \mathcal{H}^m$ such that $\Phi(T^m, h^m) \rightarrow \Phi(T, h) \,$ for $m\rightarrow \infty$. The problem is rather that the two networks compete, and thus not just their absolute approximation capabilities are relevant, but also their approximation capabilities relative to each other.
		\end{itemize}

	\end{remark}

	\section{Reformulations}
	\label{sec:reformulations}
	This section studies theoretical reformulations of problem $(P)$. The reformulations aim at improving theoretical and numerical aspects of the problem, while introducing only small changes to the objective. Among others, we show that the reformulations are better suited for approximation by neural networks and that certain aspects of the optimization are made easier by going from linear to strictly convex structures.

	First, we give two reformulations from the optimal transport literature that can loosely be described as relaxing the marginal constraints. We subsequently show how to generalize these reformulations to arbitrary problems $(P)$.
	
	\subsection{Relaxation of constraints: The optimal transport case}
	Throughout this subsection, let $\mathbb{R}^d = \mathbb{R}^{d_1} \times \mathbb{R}^{d_2}$, fix $\mu \in \mathcal{P}(\mathbb{R}^d)$ and let $\mathcal{H} := \{ h \in C_b(\mathbb{R}^{d}) : \exists \, h_1 \in C_b(\mathbb{R}^{d_1}), h_2 \in C_b(\mathbb{R}^{d_2}) : \forall (x_1, x_2) \in \mathbb{R}^{d} : h(x_1, x_2) = h_1(x_1) + h_2(x_2)\}$. This leads to $\mathcal{Q} = \Pi(\mu_1, \mu_2)$. Let $S_1$ and $S_2$ denote the projections of $\mathbb{R}^d$ onto $\mathbb{R}^{d_1}$ and $\mathbb{R}^{d_2}$, respectively.
	
	\paragraph{\cite{xie2019scalable}.} The idea of this paper is to reformulate the constraint $\mu_j = \nu_j$ as $W_1(\nu_j, \mu_j) = 0$, where $W_1(\nu, \mu)$ denotes the 1-Wasserstein distance between $\nu$ and $\mu$, and then relax this constraint with a fixed but large Lagrange multiplier $\eta > 0$. The relaxed form of $(P)$ is then given by
	\begin{equation*}
		\sup_{\nu \in \mathcal{P}(\mathbb{R}^d)}\int f \, d\nu - \eta \left( W_1\left(\nu_1, \mu_1 \right) + W_1\left(\nu_2, \mu_2  \right)   \right), 
	\end{equation*}
	for some constant $\eta > 0$. The MinMax form reduces to
	\begin{equation} \label{minimax_Xetal}
		\begin{split}
			(OT)_1 = \; & \sup_{T: \mathbb{R}^K \rightarrow \mathbb{R}^d}  \inf_{h_j \in \textrm{Lip}_1(\mathbb{R}^{d_j})} \int f\, d\theta_T - \eta  \sum_{j = 1}^{2} \left(\int h_j \circ S_j \, d\theta_T - \int h_j \, d\mu_j\right),
		\end{split}
	\end{equation} 
	with $\textrm{Lip}_1(\mathbb{R}^{d_j})$ being the 1-Lipschitz functions mapping from $\mathbb{R}^{d_j}$ to $\mathbb{R}$ for $j=1, 2$.
	
	\paragraph{\cite{yang2018scalable}.} In this paper, the constraints $\mu_j = \nu_j$ are instead penalized within the optimization problem by a $\psi$-divergence $D_\psi$. Problem $(P)$ is reformulated as
	\begin{equation*}
		\sup_{\nu \in \mathcal{P}(\mathbb{R}^d)} \int f d\nu - D_{\psi_1}(\nu_1, \mu_1) - D_{\psi_2}(\nu_2, \mu_2).
	\end{equation*}
	Utilizing the dual representations of the divergences, the MinMax form can be stated as
	\begin{equation}
		\label{eq:yanguhler}
		(OT)_2 = \sup_{T: \mathbb{R}^K \rightarrow \mathbb{R}^d} \inf_{h_j\in C(\mathbb{R}^{d_j})} \int f\,d\theta_T - \sum_{j=1}^2 \left( \int h_j \circ S_j \,d\theta_T - \int \psi_j^{\ast}(h_j) \,d\mu_j \right),
	\end{equation}
	where $\psi_j^*$ denotes the convex conjugate of $\psi_j$ for $j=1, 2$.
The problem $(OT)_2$ is an unbalanced OT problem (see, e.g., 
\cite{Chizat2018UnbalancedOT} for an overview), which 
also enables transportation between marginals which are not necessarily 
normalized to have the same mass. In the discrete case, unbalanced OT has 
computational benefits compared to the standard OT problem
(\cite{pham2020unbalanced}).
	
	\subsection{Relaxation of constraints: The general case}
	\paragraph{Lipschitz regularization.} This paragraph generalizes the regularization technique from \cite{xie2019scalable}. Denote by $\textrm{Lip}_L(\mathbb{R}^{d})$ the centered $L$-Lipschitz functions mapping from $\mathbb{R}^{d}$ to $\mathbb{R}$.\footnote{We call $f$ \textsl{centered} if $f(0)=0$. This assumption is made in an attempt to avoid trivial scaling issues later on. Notably, in the dual formulation of the Wasserstein distance, restricting to centered functions can be done without loss of generality.}
	We define $\mathcal{H}_L$ analogously to $\mathcal{H}$, except that $C(\mathbb{R}^{d_j})$ is replaced by $\textrm{Lip}_L(\mathbb{R}^{d_j})$, i.e., we set  $\mathcal{H}_L := \left\{ \sum_{j=1}^J e_j \cdot (h_j \circ \pi_j) : h_j \in \textrm{Lip}_L(\mathbb{R}^{d_j})\right\}$. Correspondingly, we define $\textrm{Lip}^m_L(\mathbb{R}^{d}) := N^m_{d,1} \cap \textrm{Lip}_L(\mathbb{R}^{d})$ and $\mathcal{H}_L^m := \left\{ \sum_{j=1}^J e_j \cdot (h_j \circ \pi_j) : h_j \in \textrm{Lip}_L^m(\mathbb{R}^{d_j})\right\}$. Notably, the set $\textrm{Lip}_L^m(\mathbb{R}^d)$ still satisfies universal approximation properties (see \cite{eckstein2020lipschitz}). Define 
	\begin{align}
	\label{eq:PLdef}
	(P_L) := \sup_{T : \mathbb{R}^{K} \to \mathbb{R}^{d}} \inf_{h \in \mathcal{H}_L} \int f\, d\theta_T + \int h \, d\theta_T - \int h\, d\mu,\\
	(P_L^m) := \sup_{T \in N_{K, d}^m} \inf_{h \in \mathcal{H}_L^m} \int f \, d\theta_T + \int h \, d\theta_T - \int h\, d\mu.
	\end{align}
	We note that for $L = \eta$, in the optimal transport case, $(P_L) = (OT)_1$.
	
	\paragraph{Divergence regularization.} This paragraph generalizes the regularization technique from \cite{yang2018scalable}. The standard MinMax formulation of $(P)$, as derived in \eqref{MM}, can be rewritten as
	\begin{equation}
		(P) = \sup_{T: \mathbb{R}^K \rightarrow \mathbb{R}^d} \inf_{h_j \in C_b(\mathbb{R}^{d_j})} \int f \,d\theta_T - \sum_{j=1}^J \Bigg( \int e_j \cdot (h_j \circ \pi_j) \,d\theta_T - \int e_j \cdot (h_j \circ \pi_j) \,d\mu \Bigg)  .
	\end{equation}
	Introducing convex functions $\psi_j : \mathbb{R} \rightarrow \mathbb{R}$ for $j=1,\dots, J$, we define
	\begin{equation}
		\label{eq:relax}
		\begin{split}
		(P_{\psi}) = \sup_{T: \mathbb{R}^K \rightarrow \mathbb{R}^d}  \inf_{h_j \in C_b(\mathbb{R}^{d_j})} \int f  \,d\theta_T\, - \sum_{j=1}^J \Bigg( \int e_j \cdot (h_j \circ \pi_j) \,d\theta_T\\ - \int e_j \cdot (h_j \circ \pi_j) + |e_j| \cdot \psi_j^{\ast}(h_j \circ \pi_j) \,d\mu \Bigg).
		\end{split}
	\end{equation}
	Analogously, we define $(P_{\psi}^m)$, with $T \in N_{K, d}^m$ and $h_j \in N^m_{d_j,1}$.
	
	We note that, in the optimal transport case, one can recover the formulation in \eqref{eq:yanguhler} by \cite{yang2018scalable} as follows: Consider in \eqref{eq:yanguhler} the divergences with convex functions $\tilde{\psi}_1$ and $\tilde{\psi}_2$. This is recovered in \eqref{eq:relax} by setting $\psi^{\ast}_1(x) = \tilde{\psi}^\ast_1(x) - x$ and $\psi^{\ast}_2(x) = \tilde{\psi}^\ast_2(x) - x$.

	We now state the main theorem which showcases approximation capabilities of neural networks for the problems $(P_L)$ and $(P_\psi)$.
	
	\begin{theorem}\label{thm:approxreg} Assume that all measures in $\mathcal{Q} \neq \emptyset$ are compactly supported on $K\subseteq \mathbb{R}^d$, $e_1,\dots, e_J, \pi_0, \dots, \pi_J$ are Lipschitz continuous and all maps $T$ are restricted to have range $K$.\footnote{Restricting $T$ to have range $K$ is understood in the sense that the actual output $\tilde{T}(x)$ of the network will be projected onto the set $K$, i.e., $T(x) := \arg\min_{y\in K} |y-\tilde{T}(x)|$. For the statement of the theorem, the only important consequence is that since $K$ is assumed to be compact, $T$ is compact-valued as well.} Assume that the activation function of the networks for $h_j$ is either one-time continuously differentiable and not polynomial, or the ReLU function.
		\begin{itemize}
			\item[(i)]  It holds $(P_L^m) \rightarrow (P_L)$ for $m \rightarrow \infty$.
			\item[(ii)] Assume $e_j \geq 0$ for $j=1, \dots, J$, $\hat{\nu} = {\theta_{\hat{T}}} \in \mathcal{Q}$ is an optimizer of $(P_\psi)$ and there exists a sequence of network functions $T_m \in N_{K, d}^m$ such that $\frac{d {\theta_{T_{m}}}}{d \theta_{\hat{T}}}$ is bounded and converges almost surely to 1. Then $\liminf_{m\rightarrow \infty}(P_\psi^m) \geq (P_\psi)$.
		\end{itemize}
	\end{theorem}
	
	\begin{remark}
		\label{rem:reg}
		Theorem \ref{thm:approxreg} showcases that the problems $(P_L)$ and 
		$(P_\psi)$ are advantageous compared to $(P)$ in the sense that the 
		neural network approximations $(P_L^m) \approx (P_L)$ and $(P_\psi^m) 
		\approx (P_\psi)$ are more justified compared to $(P^m) \approx (P)$. 
		In particular, the obstacle that the approximation is severely uneven, 
		in the sense that the inner infimum always evaluates to minus infinity 
		(as illustrated in Remark \ref{rem:pmapprox}), is remedied in both 
		situations. Recall that for $(P)$ and $(P^m)$ any violated constraint 
		can be punished arbitrarily. On the other hand, for the regularizations:
		\begin{itemize}
			\item[(i)] For $(P_L)$, for a fixed $T$, the inner infimum in Eq. \eqref{eq:PLdef} intuitively calculates a scaled (by the factor $L$) Wasserstein-like distance between $\theta_T$ and the set of feasible solutions. This means that, compared to problem $(P)$, it is now quantified \textbf{how much} a constraint is violated, and not just whether or not it is. 
			\item[(ii)] For $(P_\psi)$, a similar logic as in $(i)$ applies. In this case, however, instead of the Wasserstein distance, the inner infimum calculates a kind of divergence. Again, the deviation between $\theta_T$ and the set of feasible solutions is quantified. Theoretically, as divergences can evaluate to infinity (in particular whenever absolute continuity issues occur), this leads to worse approximation properties than when using the Wasserstein distance. However, a different advantage is given by the fact that for $(P_\psi)$ the inner infimum is now strictly convex in the function $h$, which can greatly improve stability in the numerics (see Section \ref{ex:distr_constraint}) .
		\end{itemize}
	\end{remark}
	
		\begin{remark}
			The result stated in Theorem \ref{thm:approxreg} (ii) leaves several questions open:
			\begin{itemize}
				\item[(i)] Existence of $T_m$ satisfying $\frac{d {\theta_{T_{m}}}}{d \theta_{\hat{T}}} \rightarrow 1$ is difficult to verify in general. In the case where $\theta_{\hat{T}}$ has a continuous Lebesgue density $g$ which is bounded away from zero (on the compact set $K$ considered in Theorem \ref{thm:approxreg}), the assumption can be simplified. It then suffices to assume that $T_m$ converges pointwise to $T$ (this yields $\theta_{T_m} \stackrel{w}{\rightarrow} \theta_T$) and that $\theta_{T_m}$ ($m \in \mathbb{N}$) have Lebesgue densities $g_m$ that are equicontinuous and uniformly bounded, since then $\frac{d {\theta_{T_{m}}}}{d \theta_{\hat{T}}} = \frac{g_n}{g} \rightarrow 1$ as $g_n \rightarrow g$ uniformly by \cite{boos1985converse} and using that $g$ is bounded away from zero.
				\item[(ii)] We expect that the converse, i.e., $\limsup_{m\rightarrow \infty} (P_{\psi}^m) \leq (P_\psi)$, may be shown as well for certain divergences $\psi$. The difficulty hereby is the following: For $(P_L)$, the reason the converse works is that the set of normalized Lipschitz functions is a bounded and equicontinuous set, and hence compact by the Arzel\`{a}-Ascoli theorem. The functions $h_j$ occurring in the inner infimum of $(P_\psi)$ do not satisfy such a compactness property a priori. Nevertheless, we still expect that one might effectively reduce to such a compact case for $(P_\psi)$ as well (given sufficiently nice $\psi$). A thorough analysis is however left for future work.
			\end{itemize}
		\end{remark}

	\section{Algorithmic considerations}
	\label{sec:algorithmic}
	The most widespread utilization of neural networks in MinMax settings is within GANs. This section discusses techniques from the GAN literature on overcoming instability during training when solving problem $(P)$. For discussions specific to GANs, see for instance \cite{salimans2016improved, roth2017stabilizing,thanh2019improving} and a recent survey by \cite{wiatrak2019stabilizing}.
	
	First, we mention why there is a specific need to go into detail on the training procedure for the MinMax approach for problem $(P)$, compared to just applying everything that works well within GAN settings. For GANs, the basic objective is soft, e.g., creating realistic pictures with certain features. Even if this is made rigorous (for instance via the inception score), the actual theoretical value of the MinMax problem is of little interest. Hence for GAN training one can apply procedures which change this theoretical value while improving on the other criteria of interest. Such procedures are unsuitable to apply to problem $(P)$. These include: 
	
	\begin{itemize}
		\item Batch normalization (\cite{ioffe2015batch}): With batch normalization, the respective functional spaces are altered. In the definition of $\mathcal{H}$ in \eqref{FH}, the functions $h_j$ then do not just map from $\mathbb{R}^{d_j}$, but take the batch distribution as an additional input, which can significantly change the theoretical problem.
		\item Certain forms of quantization (\cite{sinha2019small}) or noise convolution (\cite{DBLP:conf/iclr/ArjovskyB17}): Even minor adjustments to input distributions (like adding small Gaussian noise) can lead to drastic changes. Theoretically, changing the input distributions corresponds to changing the measure $\mu$ when defining problem $(P)$ in \eqref{eq:primaldef}. For certain constraints, problem $(P)$ is very sensitive to changes in the measure $\mu$: For instance, in martingale optimal transport (\cite{beiglbock2016problem}), the set $\mathcal{Q}$ may become empty with arbitrarily small changes to $\mu$.
	\end{itemize}
	
	On the other hand, some approaches that work well within the literature on GANs are certainly also feasible for problem $(P)$. These include:
	\begin{itemize} 
		\item Multi-agent GANs (\cite{ghosh2018multi,hoang2018mgan,ahmetouglu2019hierarchical}): This approach involves the introduction of multiple generators and/or discriminators. In game theoretical terms, this makes it easier for the players (generator and discriminator) to utilize mixed strategies, which is essential for obtaining stable equilibria. When utilizing mixtures for the discriminator, slight care has to be taken, since a mixture between multiple discriminators is usually not continuous, which is at odds with the space of functions used to define $\mathcal{H}$ in \eqref{FH}. Usually, however, the statement of problem $(P)$ is robust with respect to such changes. Further, mixtures for the generator can always be utilized.
		\item Variations on GDA methods: Among others, methods like \textit{unrolled} GANs (\cite{Metz2017unrolled}), \textit{consensus optimization} (\cite{mescheder2017numerics}), \textit{competitive gradient descent} (\cite{schafer2019competitive}), or the \textit{follow-the-ridge} approach (\cite{Wang2019ridge}), change the way that parameters are updated during training, and hence do not affect the theoretical objective at all, while improving the training procedure.
		\item Scaling up: As pursued in \cite{brock2018large} for GANs, solving problem $(P)$ can be improved by scaling up the size of the networks and increasing the batch size. Increasing the batch size is particularly suitable for problem $(P)$, as the ``true distribution'' $\mu$ is known and one can produce arbitrary amounts of samples.
	\end{itemize}
	A further point to be mentioned which can be useful to consider is the choice of the latent measure $\theta$. In GANs, this is usually simply taken as a high dimensional standard normal distribution. On the other hand, settings related to autoencoders and optimal transport examine the choice of the latent measure in more depth (\cite{rubenstein2018latent, henry2019martingale}).

	For problem $(P)$, taking mixtures of generators (\cite{ghosh2018multi}) and 5 to 10 unrolling steps (\cite{Metz2017unrolled}) has proven to be very well suited. In the following, we quickly argue why this is the case, while Section \ref{ex:MOT} supports this claim with numerical experiments. To this end, recall that the fundamental goal of the generator is to be able to generate a wide spectrum of possible candidate solutions. For optimal transport, it is often sufficient to generate Monge couplings (see \cite{gangbo1996geometry}), which are relatively simple. For general problems $(P)$, the required candidates often have to be more complex (see, e.g., Section 3.8 in \cite{henry2019martingale}). While a single trained generator may be biased towards concentrated measures, multiple generators more easily represent smooth measures as well. On the other hand, the fundamental goal of the discriminator is to punish the generator if the proposed candidate does not satisfy the constraints. During training, the optimization will usually reach points where the generator tries to push for a slight violation of the constraints that cannot be immediately punished by the discriminator, while slightly improving the objective value. We found that using unrolling (in which the generator takes possible future adjustments of the discriminator into account) greatly restricts violations of the constraints, because the generator already anticipates punishment in the future even at the current update of its parameters.

\section{Numerical experiments}
\label{sec:numerics}
Code to reproduce the numerical experiments is available on \url{https://github.com/stephaneckstein/minmaxot}.
\subsection{Optimal transport with distribution constraint (DCOT)} \label{ex:distr_constraint}

In this section, we consider an optimal transport problem with an additional distribution constraint. Let $f(x_1, x_2) := (x_1+x_2)^+$, $\mu_1, \mu_2$ be normal distributions with mean 0 and standard deviation 2, and $\kappa \in \mathcal{P}(\mathbb{R})$ a Student's $t$-distribution with 8 degrees of freedom. Set $\mathcal{H} = \{h_1(x_1) + h_2(x_2) + u(x_2-x_1) : h_1, h_2, u \in C_b(\mathbb{R})\}$. We use the notation $\bar\pi(x_1, x_2) = (x_2-x_1)$. Choosing $\mu$ with the specified marginals such that $\bar{\pi}_*\mu = \kappa$, this leads to ${\mathcal{Q} = \left\{ \nu \in \mathcal{P}(\mathbb{R}^2) : \nu_1 = \mu_1 \text{ and } \nu_2 = \mu_2 \mbox{\, such that \,} \bar{\pi}_*\nu = \kappa \right\}}$.\footnote{The described problem may occur naturally in financial contexts, where two assets have described distribution $\mu_1$ and $\mu_2$. The function $f$ models the payoff of a basket option and the constraint including $\kappa$ may describe information about the relation between the two assets.}

For additional details on the specifications, see Appendix \ref{app:ex_dist}.

\begin{figure}[t] 
	\centering 
	\includegraphics[width=\textwidth]{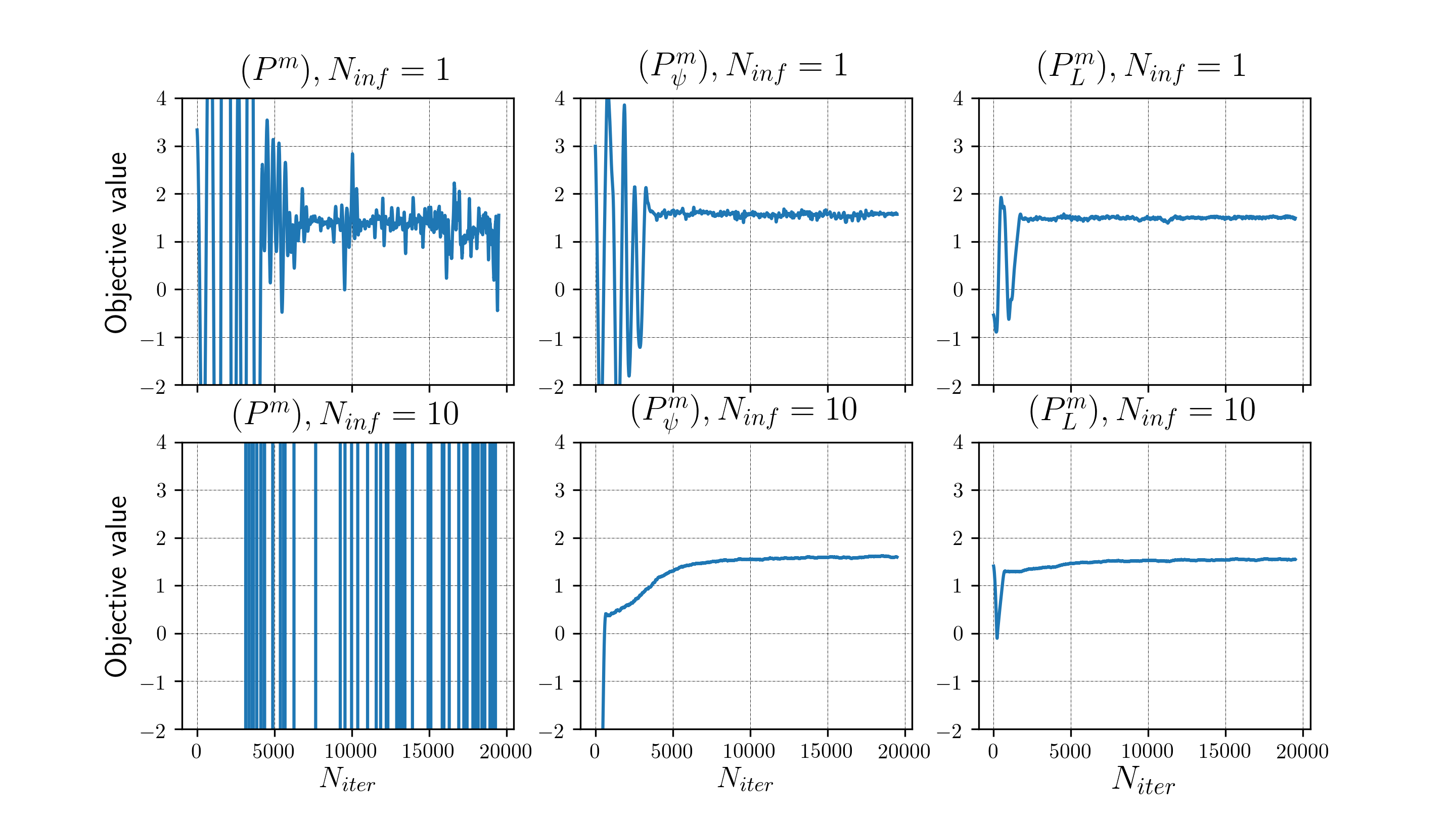} 
	\caption{Results for Section \ref{ex:distr_constraint}. Numerical 
	convergence observed for the optimal transport problem with additional 
	distribution constraint as described in Section \ref{ex:distr_constraint} 
	and Appendix \ref{app:ex_dist}. The top row shows the convergence when 
	parameters are updated using a GDA algorithm with a single updating step 
	for both infimum and supremum. The bottom row shows the convergence in the 
	case of 10 infimum updates for each supremum update. The first column 
	reports the case where no regularization is used (problem $(P^m)$). The 
	second and third column report the results for the regularization methods 
	$(P_\psi)$ and $(P_L)$, respectively, as described in \mbox{Section 
	\ref{sec:reformulations}}. The displayed graphs are median values across 11 
	runs with respect to the standard deviation of the objective values over 
	the last 5000 iterations.} \label{fig:distr_constraint}
\end{figure}
This experiment is meant to showcase the insights obtained in Remark \ref{rem:pmapprox} and Theorem \ref{thm:approxreg}. We first calculate the problems $(P^m), (P_\psi^m)$ and $(P_L^m)$ using GDA with Adam optimizer, with a single alternating updating step for both supremum and infimum networks. We observe the convergence and stability in the top row of Figure \ref{fig:distr_constraint}. While the graphs including regularization showcase better stability, the real benefit is revealed in the bottom row of Figure \ref{fig:distr_constraint}. When increasing to 10 the number of infimum steps in the GDA, the implemented problem more closely resembles the theoretical one, since now the infimum can really be regarded as the inner problem. As predicted by Remark \ref{rem:pmapprox}, the calculation of $(P^m)$ is now entirely unstable. On the other hand, consistent with Remark \ref{rem:reg}, the convergence for the problems $(P^m_\psi)$ and $(P^m_L)$ becomes even more smooth.

\subsection{Martingale optimal transport (MOT)}
In this section, the martingale optimal transport problem (\cite{beiglbock2013model}) is studied.  We consider the cost function $f(x_1, x_2) = (x_2 -x_1)^{+}$. Also, let $\mathcal{N}(m, \sigma^2)$ describe a normal distribution with mean $m$ and variance $\sigma^2$. We define the marginals as follows:
\begin{equation*} 
\begin{split}
& \mu_1 = 0.5 \, \mathcal{N}(-1.3, 0.5^2) + 0.5 \, \mathcal{N}(0.8, 0.7^2), \\
& \mu_2 = 0.5 \, \mathcal{N}(-1.3, 1.1^2) + 0.5 \, \mathcal{N}(0.8, 1.3^2).
\end{split}
\end{equation*}
Set $\mathcal{H} = \{h_1(x_1) + h_2(x_2) + g(x_1)\cdot(x_2-x_1) : h_1, h_2, g \in C_b(\mathbb{R})\}$. Then, we obtain $ {\mathcal{Q} = \{ \nu \in \mathcal{P}(\mathbb{R}^2): \nu_1 = \mu_1, \nu_2 = \mu_2 \text { and, if } (X_1, X_2) \sim \nu, \text{ then } \mathbb{E}[X_2 | X_1] = X_1\}}$.\footnote{In financial terms, this example corresponds to computing price bounds on a forward start call option under a martingale constraint. For related problems, see, e.g., \cite{beiglbock2013model}.}

For additional details on the specifications, see Appendix \ref{app:ex_mot}.

\label{ex:MOT}
In the first row of Table \ref{table:exmot}, we observe how a standard alternating updating of generator and discriminator parameters can lead to difficulties with respect to both stability of the convergence and feasibility of the obtained solution. We then resolve these issues by adjusting the algorithmic procedure and report the results in the bottom rows of Table \ref{table:exmot}. These results show that using a mixture of generators or utilizing unrolling can greatly improve stability and feasibility issues, as well as improve the optimal value of the obtained solution. Combining the two methods leads to the best results.
\begin{table}
	\caption{Results for Section \ref{ex:MOT}}
	\label{table:exmot}
	\begingroup
	\centering
	\begin{tabular}{l  l l l l}
		\toprule
		& Integral value& Error marginals & Error martingale & Std dev iterations\\\midrule
		Base & 0.281 & 0.126 & 0.087 & 0.267 \\
		Mixtures & 0.291 & 0.062 & 0.038 & 0.078 \\
		Unrolling & 0.289 & 0.016 & 0.011 & 0.025 \\
		Combined & \textbf{0.299} & \textbf{0.014} & \textbf{0.010} & \textbf{0.015} \vspace{0.2cm} \\\bottomrule
	\end{tabular}\\\vspace{1.5mm}
	\endgroup
	\noindent
	Average values obtained over 10 runs of solving the MOT problem as 
	described in Section \ref{ex:MOT} and Appendix \ref{app:ex_mot}. The 
	unrolling procedure is performed with 5 unrolling steps of the 
	discriminator, and for the mixing, a fixed mixture of 5 generators is used. 
	The ``Combined'' row uses both mixtures and unrolling. As the problem is a 
	maximization problem, high integral values while having low error values 
	are desirable. The error values hereby quantify violations of the 
	constraints. The final column gives an indicator for the stability during 
	training, where low values imply stable convergence.	
\end{table}

\section{Conclusion and outlook}
We introduced a general MinMax setting for the class of problems of the form $(P)$. We argued that regularization techniques known from the OT literature can be generalized. By proving approximation theorems, we gave theoretical justification for utilizing neural nets to calculate the solution of the regularized problems. Further, we argued that, with regularization, the inner infimum of the MinMax problem is usually bounded and thus instability during training can be reduced. Beyond the theoretical objective, we discussed algorithmic adjustments that can be adapted from the GAN literature. Both theoretical insights - utilizing regularization and algorithmic adjustments - were shown to be beneficial when applied in numerical experiments.

The following avenues for future research are left open. Firstly, some aspects of the theoretical approximations introduced by $(P_L)$ and $(P_\psi)$ can be studied in more depth (in particular, a rigorous analysis on the approximation errors $\vert (P_L) - (P) \vert $ and $\vert (P_\psi) - (P) \vert $, see also Appendix \ref{app:approx}). Secondly, a thorough comparison with existing methods on large scale problems can give further insights on the computational possibilities. And finally, quantitative rates of the convergences studied in Theorem \ref{thm:approxreg} are of practical interest.

\section*{Acknowledgements}
The authors like to thank Michael Kupper, Carole Bernard and the referees for stimulating discussions and helpful remarks.
Luca De Gennaro Aquino is grateful to Michael Kupper and Stephan Eckstein for their hospitality at the University of Konstanz, where part of this project was done. Stephan Eckstein is thankful to the Landesgraduiertenförderung Baden-Württemberg for financial support.

\newpage
\section*{Broader Impact}
In this paper, we formally provide justification for utilizing neural networks when solving a frequently used class of optimization problems. We believe that our results can function as theoretical and practical guidelines for researchers (and practitioners) who are interested in exploring possible applications of optimal transport and related frameworks utilizing MinMax methods. 

However, it is important to emphasize that, generally speaking, theoretical insights might still be restricted by numerical convergence, thus we do not encourage overconfidence in the solution methods when resorting to neural networks.  

Nonetheless, we do not expect our work to feasibly induce any disadvantage for any group of people, nor that particular consequences for the failure of the proposed optimization methods might occur.

\bibliography{refs}
\bibliographystyle{abbrvnat}

\newpage
\appendix
\section{Proofs} \label{Appendix:proofs}
\paragraph{Proof of Theorem \ref{thm:approxreg}.} Throughout, we use the notation $\Phi(T, h) := \int f \circ T \, d\theta + \int h \circ T  \, d\theta - \int h\, d\mu = \int f \,d\theta_T + \sum_{j=1}^J \int e_j \cdot (h_j \circ \pi_j) \,d\theta_T - \int e_j \cdot (h_j \circ \pi_j) \,d\mu$. 

Proof of (i): We show that, for a given $\varepsilon > 0$, there is $m\in \mathbb{N}$ such that both $(P_L) \stackrel{(a)}{\geq} (P_L^m) - \varepsilon$ and $(P_L) \stackrel{(b)}{\leq} (P_L^m) - \varepsilon$ hold. 

Regarding (a), choose $m$ such that any $L$-Lipschitz function on the compact set $K$ can be approximated up to accuracy $\varepsilon/(2J\max_{j=1,\dots, J} \|e_j\|_\infty)$ in $\|\cdot\|_\infty$ by neural networks with hidden dimension $m$, which is possible by \cite[Theorem 1]{eckstein2020lipschitz}. Then, for all $T$, it holds for any $j = 1,\dots, J$ that
\[
\inf_{h_j \in \textrm{Lip}_L^m} \int e_j \cdot (h_j \circ \pi_j) \,d\theta_T + \int e_j \cdot (h_j \circ \pi_j) \,d\mu - \inf_{h_j \in \textrm{Lip}_L} \int e_j \cdot (h_j \circ \pi_j) \,d\theta_T + \int e_j \cdot (h_j \circ \pi_j) \,d\mu \leq \varepsilon/J
\]
and thus
$\inf_{h\in \mathcal{H}_L} \Phi(T, h) \geq \inf_{h\in \mathcal{H}_L^m} \Phi(T, h) - \varepsilon$. This implies \[(P_L) \geq \sup_{T \in N_{K, d}^m} \inf_{h \in \mathcal{H}_L} \Phi(T, h) \geq \sup_{T \in N_{K, d}^m} \inf_{h \in \mathcal{H}_L^m} \Phi(T, h) - \varepsilon = (P_L^m) - \varepsilon,\] and hence (a) follows. 

Regarding (b), choose an optimizer $\hat{\nu} = \hat{T}_*\theta$ of $(P_L)$. Since $T$ is compact-valued, $T\in L_1(\theta)$, and hence we can choose $T_m \in N_{K, d}^m$ such that $T_m \rightarrow T$ for $m\rightarrow \infty$ in $L_1(\theta)$, which implies $\nu^m := (T_m)_*\theta \stackrel{w}{\rightarrow }\hat{\nu}$ and since the measures are supported on $K$ also $W_1(\nu^m, \hat{\nu}) \rightarrow 0$ for $m\rightarrow \infty$. It holds
\[\bigg|\inf_{h \in \mathcal{H}_L} \phi(\hat{T}, h) - \inf_{h \in \mathcal{H}_L} \phi(T^m, h)\bigg| \leq \sum_{j=1}^J \sup_{h_j \in \rm{Lip}_L(\mathbb{R}^{d_j})} \Bigg| \int e_j \cdot (h_j \circ \pi_j) \,d\hat{\nu} - \int e_j \cdot (h_j \circ \pi_j) \,d\nu^m \Bigg| =: (*)\]

Note further that there exists some $\hat{L} > 0$ such that all $h_j \circ \pi_j$ are $\hat{L}$-Lipschitz. Since $h_j$ are centered and compact-valued, their infinity norms are bounded uniformly, say by some $S>0$. Hence any $e_j \cdot (h_j \circ \pi_j)$ is $(\hat{L} \|e_j\|_\infty + L_{e_j} S)$-Lipschitz. We denote the maximum of these constants by $\bar{L}$. Thus $(*) \leq J \bar{L} W_1(\nu^m, \hat{\nu}) \leq \varepsilon/2$ for $m$ large enough. Also, $|\int f \,d\hat{\nu} - \int f \,d\nu^m|\leq \varepsilon/2$ for $m$ large enough, since $f$ restricted to $K$ is continuous and bounded. Hence
\[
(P_L) = \inf_{h \in \mathcal{H}_L} \Phi(\hat{T}, h) \leq \inf_{h \in \mathcal{H}_L} \Phi(T^m, h) + \varepsilon \leq \inf_{h \in \mathcal{H}_L^m} \Phi(T^m, h) + \varepsilon \leq (P_L^m) + \varepsilon,
\]
which yields the claim.

Proof of (ii): The proof builds heavily on the fact that $e_j$ are assumed to be non-negative, which allows for a reformulation of $(P_\psi)$ in terms of divergences.
For $\nu \in \mathcal{P}(\mathbb{R}^d)$, we define the measure $\nu^{e_j}$ by $\frac{d\nu^{e_j}}{d\nu} = e_j$. We get
\begin{align*}
&~~~~\sup_{h_j \in C_b(\mathbb{R}^{d_j})} \int e_j \cdot (h_j \circ \pi_j) \,d\theta_T - \int \left( e_j \cdot (h_j \circ \pi_j) + e_j \cdot \psi_j^{\ast}(h_j \circ \pi_j) \right) \,d\mu \\
&= \sup_{h_j \in C_b(\mathbb{R}^{d_j})} \int h_j \circ \pi_j \,d\theta_T^{e_j} - \int \left( h_j \circ \pi_j + \psi_j^{\ast}(h_j \circ \pi_j) \right) \,d\mu^{e_j} \\
&= \sup_{h_j \in C_b(\mathbb{R}^{d_j})} \int h_j \,d\big((\pi_j)_*\theta_T^{e_j}\big) - \int \left( h_j + \psi_j^{\ast}(h_j) \right) \,d\big((\pi_j)_*\mu^{e_j}\big) \\
&= D_{\tilde{\psi}_j}\big((\pi_j)_*\theta_T^{e_j}, (\pi_j)_*\mu^{e_j}\big),
\end{align*}
where $\tilde{\psi}_j^\ast(x) = x + \psi_j^\ast(x)$, $D_{\tilde{\psi}_j}(\nu, \mu) = \int \tilde{\psi}\big(\frac{d\nu}{d\mu}\big) d\mu$ for $\nu \ll \mu$, and the last equality follows by the dual representation for divergences.\footnote{See for instance \cite[Chapter 4]{broniatowski2006minimization}, and note that while the dual formulation therein is based on bounded and measurable functions, on the compact set $K$ standard approximation arguments using Lusin's and Tietze's theorems yield that continuous functions are sufficient.}
The above shows that 
\[
(P_{\psi}) = \sup_{T : \mathbb{R}^K \rightarrow \mathbb{R}^d} \int f \,d\theta_T - \sum_{j=1}^J D_{\tilde{\psi}_j}\big((\pi_j)_*\theta_T^{e_j}, (\pi_j)_*\mu^{e_j}\big).
\]
Now, choose an optimizer $T$ and a sequence $T^m \in N_{K, d}^m$ as in the 
assumption of the theorem. Without loss of generality, we can choose a 
representative among the almost-sure equivalence class, such that 
$\frac{d\theta_{T^m}}{d\theta_T} \rightarrow 1$ holds point-wise for 
$m\rightarrow \infty$. Elementary calculation yields that 
$\frac{(\pi_j)_*\theta_{T^m}^{e_j}}{(\pi_j)_*\theta_T^{e_j}} \rightarrow 
1$ holds point-wise as well, and hence by dominated convergence 
$D_{\tilde{\psi}_j}((\pi_j)_*\theta_{T^m}^{e_j}, (\pi_j)_*\mu^{e_j}) \rightarrow D_{\tilde{\psi}_j}((\pi_j)_*\theta_T^{e_j}, 
(\pi_j)_*\mu^{e_j})$ for $m\rightarrow \infty$ follows. We can choose $m 
\in \mathbb{N}$ such that $(P_\psi) \leq \int f \,d\theta_{T^m} - \sum_{j=1}^J 
D_{\tilde{\psi}_j}((\pi_j)_*\theta_{T^m}^{e_j}, (\pi_j)_*\mu^{e_j}) + \varepsilon$. By again plugging in the dual formulation for 
$D_{\tilde{\psi}_j}$, and noting that the infimum only gets larger when 
restricted to neural network functions,
\begin{align*}
(P_{\psi}) \leq & \inf_{h_j \in N_{d_j, 1}^m} \int f \,d\theta_{T^m}  -  \sum_{j=1}^J \Bigg( \int e_j \cdot (h_j \circ \pi_j) \,d\theta_{T^m} \\ 
& - \int e_j \cdot (h_j \circ \pi_j) + e_j \cdot \psi_j^{\ast}(h_j \circ \pi_j) \,d\mu \Bigg) + \varepsilon \leq (P_{\psi}^m) + \varepsilon,
\end{align*}
which yields the claim.

\section{Specifications of numerical examples}
\label{Appendix:Examples}
Here we provide a quick overview of the specifications for the numerical experiments discussed in Section \ref{sec:numerics}. Further details can be seen within the code on \url{https://github.com/stephaneckstein/minmaxot}. 

In all examples, we use the Adam optimizer (\cite{kingma2014adam}) with learning rate ${\alpha = 10^{-5}}$ and $\beta_1 = 0.5, \,\beta_2 = 0.999$ and $\epsilon = 10^{-9}$. 
Both generator and discriminator consist of 4 layer feed-forward networks with hidden dimension 64 (for Section \ref{ex:distr_constraint}) or 128 (for Section \ref{ex:MOT}). Network weights are initialized using the GlorotNormal initializer. For the generator networks, we choose the hyperbolic tangent activation function. For the discriminator networks, we choose the ReLU activation function. Computations are performed in Python 3.7 using TensorFlow 1.15.0.

\subsection{Specification of the experiment in Section \ref{ex:distr_constraint}} \label{app:ex_dist}

For $(P_L)$, we take $L=1$, and implement the Lipschitz constraint as described 
in Appendix \ref{subsec:gradientpenalty}. Although $L=1$ appears low, since $f$ 
is also 1-Lipschitz we found this choice to be sufficient. If $L$ is chosen 
larger, the obtained objective value does not appear to change significantly, 
but the stability during training gets slightly worse. For $(P_\psi)$, we take 
$\psi_j^*(x) = \frac{x^2}{25}$ for $j=1, 2, 3$, and we found other choices 
(see, e.g., Table 1 of \cite{yang2018scalable} for a list of candidates) to be 
comparable regarding the improved stability during training. For intuition 
regarding both choices, see also Appendix \ref{app:approx}. 

As latent measure, 
we choose $\theta = \mathcal{U}([-1, 1]^2)$ (the uniform distribution on $[-1, 
1]^2$).

The graphs in Figure \ref{fig:distr_constraint} are constructed as follows: For each supremum iteration $t$ of Algorithm \ref{alg:P}, we evaluate and save the term $\frac{1}{\min\{t, N_r\}} \sum_{s=N-\min\{t, N_r\}+1}^N \Phi_s^m(f \, ; \textbf{w}_{h}, \textbf{w}_{T})$ (where $N_r$ is set to 500), which would be the output value of the algorithm if iteration $t$ were the final iteration. The resulting list of values in dependence on the iteration is plotted in the graphs.

\subsection{Specification of the experiment in Section \ref{ex:MOT}} \label{app:ex_mot} 

As latent measure, we choose $\theta = \mathcal{U}([-1, 1]^2)$.

The first column in Table \ref{table:exmot} describes the integral value of the numerical optimizer, i.e., if $T(\cdot \, ; \textbf{w}_{T})$ is the fully trained network from Algorithm \ref{alg:P} (we chose $N=15000$, $N_r = 500$, $N_{inf} = 1$, $N_s = 0$), then the first column reports $\int f \,d\theta_{T(\cdot \, ; \textbf{w}_{T})}$ approximated using $10^5$ many samples. The second and third column are explained in Section \ref{Appendix:feasible}. 
The final column reports the standard deviation of the values $\Phi_t^m(f \, ; \textbf{w}_{h}, \textbf{w}_{T})$ for $t=N-2499, ..., N$ given within Algorithm \ref{alg:P}, which characterizes the stability of the convergence.

\begin{table}
	\caption{Runtimes for the numerical experiments}
	\begingroup
	\centering
	\parbox{.5\linewidth}{
		\centering
		\begin{tabular}{l  l }
			\toprule
			DCOT, Section \ref{ex:distr_constraint} & Runtime  \\\midrule
			$(P^{m}), N_{inf} = 1$ & 83  \\
			$(P^{m}), N_{inf} = 10$ & 464  \\
			$(P_{\psi}^{m}), N_{inf} = 1$ & 85   \\
			$(P_{\psi}^{m}),N_{inf} = 10$ & 481   \\
			$(P^{m}_{L}), N_{inf} = 1$ & 137  \\
			$(P^{m}_{L}), N_{inf} = 10$ & 909 \\ \bottomrule	
		\end{tabular}
	}
	\hfill
	\parbox{.5\linewidth}{
		\centering
		\begin{tabular}{l  l} \toprule
			MOT, Section \ref{ex:MOT} & Runtime   \\\midrule
			Base & 551  \\
			Mixtures & 918    \\
			Unrolling & 4546    \\
			Combined & 4901 \\	\bottomrule	
		\end{tabular}
	}
	\endgroup
	\vspace{1.5mm}
	
	\noindent
	Reported runtimes (in seconds) for the different experiments in Section \ref{sec:numerics}. All programs were run using an intel Core i7-6700HQ CPU@2.60GHz.
\end{table}

\subsection{Algorithm}
Algorithm \ref{alg:P} shows how to compute problem $(P^{m})$ using GDA and the Adam optimizer. The returned value yields the proxy value for $(P^m)$. The fully optimized function $T(\cdot\, ;  \textbf{w}_{T})$ serves as the approximate supremum optimizer of $(P^m)$ in the MinMax setting. Hence $\theta_{T(\cdot\, ; \textbf{w}_{T})}$ is the numerically obtained optimal measure maximizing $(P^m)$. 

The problems $(P_L^m)$ and $(P^m_\psi)$ are implemented accordingly, while only the terms $\Phi^m$ and $\Phi^m_t$ are altered. Namely, for $(P^m_\psi)$, we add the divergence terms as given in \eqref{eq:relax}, while for $(P_L^m)$, we add the gradient penalty as described in Section \ref{subsec:gradientpenalty}. To include the unrolling procedure and/or the mixture of generators, adjustments according to \cite{Metz2017unrolled} and/or \cite{ghosh2018multi} have to be included.
\begin{algorithm}[!h]
	\caption{MinMax optimization for OT and beyond: Problem $(P^m)$}
	\begin{algorithmic}
		\STATE \textbf{Inputs:} cost function $f$; measure $\mu$;  latent measure $\theta$; batch size $n$; total number of iterations $N$; number of infimum steps $N_{inf}$; number of steps for return value $N_r$; number of warm-up steps $N_s$.\\
		\textbf{Require:} random initialization of neural network weights $\textbf{w}_{h}, \textbf{w}_{T}$.
		\FOR{$t = 1,\dots, N$} \vspace{0.05cm}
		\IF{$t > N_s$}
		\FOR{$\_ = 1,\dots, N_{inf}$} \vspace{0.05cm}
		\STATE sample $\{x_i\}_{i=1}^{n} \sim \mu$ \vspace{0.05cm}
		\STATE sample $\{y_i\}_{i=1}^{n} \sim \theta$ \vspace{0.05cm}
		\STATE evaluate $ \Phi^{m}(f\, ; \textbf{w}_{h}, \textbf{w}_{T}) = \frac{1}{n}\sum_{i=1}^{n} \Big(f(T(y_i \, ;  \textbf{w}_{T})) \, +  h(T(y_i \, ; \textbf{w}_{T})\, ; \textbf{w}_{h}) \,  -  h(x_i \, ; \textbf{w}_{h})\Big)$ 
		\STATE  $\textbf{w}_{h} \leftarrow$ $\mbox{Adam}(\Phi^{m}(f\, ; \textbf{w}_{h}, \textbf{w}_{T}))$ \vspace{0.05cm}
		\ENDFOR
		\ENDIF
		\STATE sample $\{x_i\}_{i=1}^{n} \sim \mu$ \vspace{0.05cm}
		\STATE sample $\{y_i\}_{i=1}^{n} \sim \theta$ \vspace{0.05cm}
		\STATE evaluate $ \Phi_t^{m}(f\, ; \textbf{w}_{h}, \textbf{w}_{T}) = \frac{1}{n}\sum_{i=1}^{n} \Big(f(T(y_i \, ; \textbf{w}_{T})) \, +  h(T(y_i\, ; \textbf{w}_{T})\, ; \textbf{w}_{h}) \,  -  h(x_i \, ; \textbf{w}_{h})\Big)$ 
		\STATE  $\textbf{w}_{T} \leftarrow$ $\mbox{Adam}(-\Phi_t^{m}(f \, ; \textbf{w}_{h}, \textbf{w}_{T}))$ \vspace{0.05cm} \vspace{0.05cm}
		\ENDFOR
		\STATE \textbf{Return:} $\frac{1}{N_r} \sum_{s=N-N_r+1}^N \Phi_s^m(f \, ; \textbf{w}_{h}, \textbf{w}_{T})$
	\end{algorithmic}
	\label{alg:P}	
\end{algorithm} 

\subsection{Numerical evaluation of feasibility}
\label{Appendix:feasible}
The numerical optimal measure $\hat{\nu} := \theta_{T(\cdot \, ; \textbf{w}_{T})}$ as given by Algorithm \ref{alg:P} should theoretically lie in $\mathcal{Q}$. To test this numerically, we (approximately) evaluate the feasibility constraint "$\, \forall h\in\mathcal{H}: \int h \,d\mu = \int h \, d\hat{\nu} \,$" for a subset of test functions $h$.

For Table \ref{table:exmot} in Section \ref{ex:MOT}, we use the first 50 Chebyshev polynomials $g_1, \dots, g_{50}$ normalized to the interval $[-6, 6]$, instead of $[-1, 1]$. For the marginal errors, the reported error is the sum $\frac{1}{2}\sum_{i=1}^2 \frac{1}{50}\sum_{j=1}^{50} |\int g_j \,d\mu_i - \int g_j \,d\hat{\nu}_i|$ where all integrals are approximated using $10^5$ many sample points. Similarly, the martingale error is the sum $\frac{1}{50}\sum_{j=1}^{50} |\int g_j(x_1)\cdot (x_2-x_1) \,\hat{\nu}(dx_1, dx_2)|$.

\subsection{Modeling Lipschitz functions}
\label{subsec:gradientpenalty}
Two methods have shown to be prevalent in the literature to enforce Lipschitz continuity: Gradient penalty (\cite{gulrajani2017improved}) and spectral normalization (\cite{miyato2018spectral}). 
We found that for our purposes a one-sided gradient penalty works well. 
To this end, enforcing $h_j \in \textrm{Lip}_L(\mathbb{R}^{d_j})$ is done via adding the penalty term
\[
\lambda \int \left((\|\nabla h_j\| - L)^+\right)^2 \,d\big((\pi_j)_*\mu\big) + \lambda \int \left((\|\nabla h_j\| - L)^+\right)^2 \,d\big((\pi_j)_*\theta_T\big),
\]
for some $\lambda > 0$, where $\|\cdot\|$ denotes the Euclidean norm.

\section{Theoretical approximations of $(P)$ by $(P_L)$ and $(P_\psi)$}
\label{app:approx}
For completeness, an analysis of the approximation of $(P)$ by $(P_L)$ and $(P_\psi)$ is required. While a full analysis is beyond the scope of this paper, we still state fundamental results:
\begin{remark}
	The definitions of $(P_L)$ and $(P_\psi)$ immediately reveal the following:
	\begin{itemize}
		\item[(i)] For $L_1 \leq L_2$, it holds $(P_{L_1}) \geq (P_{L_2}) \geq (P)$.
		\item[(ii)] For $\tilde{\psi}^*_j \geq \psi^*_j \geq 0$, $j=1, \dots, J$, it holds $(P_{\tilde{\psi}}) \geq (P_\psi) \geq (P)$.
		\item[(iii)] For $(P_L)$ to be a sensible approximation to $(P)$, $f$ has to be of linear growth, i.e., $f(x)/(1+|x|)$ has to be bounded (or even stronger restrictions have to be imposed). Otherwise it may hold $(P_L) = \infty$ for all $L$, while $(P)$ is finite. E.g., a classical OT problem on $\mathbb{R}^2$ with cost function $f(x) = |x_2-x_1|^2$ exhibits this behavior. On the other hand, numerical experiments indicate that whenever $f, \pi_1,\dots, \pi_J, e_1,\dots, e_J$ are Lipschitz continuous, it may hold $(P) = (P_L)$ for finite $L$ (see, e.g., Section \ref{ex:distr_constraint}).
	\end{itemize}
\end{remark}

\begin{table}[t!]
	\begin{adjustwidth}{-2.5cm}{-2.5cm}
		\caption{Variations of problems of the form $(P)$ from the literature}
		\label{table:listproblems}
		\begin{tabular}{l  l l}
			\toprule
			Description & $\mathcal{H}$ & Reference \\
			\midrule
			Static basket options & $\left\{ \sum_{i=1}^n \alpha_i (w_i^T x - 
			K_i)^+ : \alpha_i \in \mathbb{R}\right\}$ & \cite{d2006static} \vspace{0.1cm}\\
			Moment-constrained DRO & $\left\{c + \alpha x + \beta x^2 : c, \alpha, \beta \in \mathbb{R}\right\}$ &  \cite{popescu2007robust} \vspace{0.1cm} \\
			Optimal transport (OT) & $\left\{h_1(x) + h_2(y) : h_1, h_2 \in C_b(\mathbb{R})\right\}$ & \cite{villani2008optimal} \vspace{0.1cm}\\
			Symmetric OT & \begin{tabular}{@{}l@{}} $\Big\{\sum_{i=1}^d h_i(x_i) + (g(x_1, ..., x_d) - g(x_2, ..., x_d, x_1)): \Big. $  \\  $ \Big.  h_i \in C_b(\mathbb{R}), g \in C_b(\mathbb{R}^d) \Big\}$ \end{tabular} & \cite{ghoussoub2012remarks} \vspace{0.1cm}\\
			Martingale OT & $\left\{h_1(x) + h_2(y) + g(x)\cdot(y-x): h_1, h_2, g \in C_b(\mathbb{R}) \right\}$ & \cite{beiglbock2013model} \vspace{0.1cm}\\
			Causal OT & See Prop.~2.4 in \cite{backhoff2017causal} & \cite{lassalle2013causal}\\
			Multi-marginal OT & $\left\{\sum_{i=1}^d h_i(x_i) : h_i \in C_b(\mathbb{R})\right\}$ & \cite{pass2015multi} \vspace{0.1cm} \\
			Multi-martingale OT & \begin{tabular}{@{}l@{}} $\Big\{\sum_{i=1}^d 
				(h_{1, i}(x_i) + h_{2, i}(y_i) + g_i(x_1, ..., x_d)\cdot(y_i - x_i)) :  
				\Big. $\\ $ 
				\Big. h_{t, i} \in C_b(\mathbb{R}), g_i \in C_b(\mathbb{R}^d)\Big\}$ 
			\end{tabular} & \cite{lim2016multi} \vspace{0.1cm} \\
			OT with basket constraints & $\left\{h_1(x) + h_2(y) + c (x+y-K)^+ : 
			h_1, h_2 \in C_b(\mathbb{R}), c \in \mathbb{R}\right\}$ & 
			\cite{aquino2019bounds} \vspace{0.1cm}\\
			Finite calls MOT & \begin{tabular}{@{}l@{}}$\Big\{c + \sum_{i=1}^{n_1} \alpha_{i, 1} (x-K_{i, 1})^+ + \sum_{i=1}^{n_2} \alpha_{i, 2} (y-K_{i, 2})^+  \Big. $ \\ $ \Big. + g(x) \cdot(y-x) :c, \alpha_{i, j} \in \mathbb{R}, g \in C_b(\mathbb{R})\Big\}$ \end{tabular} & \cite[Section 3.3]{eckstein2019robust} \vspace{0.1cm} \\
			Directional OT & $\left\{h_1(x) + h_2(y) + c \eins_{\{y > x\}} : h_1, h_2 \in C_b(\mathbb{R}), c \in \mathbb{R}\right\}$ & \cite{nutz2020directional}\\
			\bottomrule \vspace{0.3cm}
		\end{tabular}
	\end{adjustwidth}
\end{table}

\section{List of problems of the form $(P)$}
\label{subsec:listofproblems}
Table \ref{table:listproblems} lists several instances of problems of the form 
$(P)$ and how they fit into the framework of this paper, i.e., how the set 
$\mathcal{H}$ is chosen. Notably, we list the simplest representatives, which 
means, for instance, in optimal transport we list the case with one dimensional 
marginal distributions. A similar class of problems as $(P)$ is studied in 
\cite{ekren2018constrained,eckstein2019computation,zaev2015monge}. 

\section{2-Wasserstein distance in $\mathbb{R}^{d}$}
\label{app:w2}
In this section, we consider the problem of computing the 2-Wasserstein distance in $\mathbb{R}^{d}$. To do so, we set the cost function $f(x_1, x_2) = - \sum_{i = 1}^{d}\left(x_{1}^{i} - x_{2}^{i}\right)^{2}.$ The marginal distributions $\mu_1$ and $\mu_2$ are chosen to be uncorrelated Gaussian distributions in $\mathbb{R}^{d}$ with means 0 and variances 1 and 4, respectively.\footnote{A similar example is discussed in  \cite{henry2019martingale}, Section 4.1.} In this case, the exact 2-Wasserstein distance is given by $W_2(\mu_1, \mu_2) = d$.  

\begin{table}[h!]
		\begin{adjustwidth}{-2.5cm}{-2.5cm}
	\caption{2-Wasserstein distance in $\mathbb{R}^{d}$}
	\label{W2}
	\begingroup
	\centering
	\begin{tabular}{l l l  l l l l}
	\toprule
	& \multicolumn{2}{c}{$(P^{m})_{base}$} & 
	 \multicolumn{2}{c}{ $(P^m)$} & 
	 \multicolumn{2}{c}{$(P^{m}_{\psi})$}  \\
	& Objective value & Std dev iterations& Objective value & Std dev iterations & Objective value & Std dev iterations \\\midrule
	& \multicolumn{6}{c}{$m = 64$} \\
	d & & & & &  &  \\
	1 & 1.055& 0.081 & 0.998 & 0.003 & 0.972 & 0.004  \\
	2 & 3.810 & 1.944 & 2.001 & 0.003 & 1.927 & 0.004 \\
	3 & 4.346 & 1.882 & 3.004 & 0.009 & 2.901 & 0.010 \\
	5 & 8.007 & 3.673 & 5.292 & 0.201 & 4.922 & 0.020 \\
	10 & 19.371 & 9.854 & 10.061 & 0.654 & 10.070 & 0.067 \vspace{0.2cm} \\ 
	& \multicolumn{6}{c}{$m = 256$} \\ 
	d & &  &  & & & \\
	1 & 1.110 & 0.285 & 1.000 & 0.003 & 1.024 & 0.007 \\
	2 & 2.048 & 0.057 & 2.004 & 0.007 & 2.026 & 0.009 \\
	3 & 3.076 & 0.093 & 2.998 & 0.006 & 3.002 & 0.012 \\
	5 & 5.359 & 0.177 & 4.993 & 0.005 & 5.028 & 0.015 \\
	10 & 16.396 & 1.800 & 9.997 & 0.008  & 10.035  & 0.015 \vspace{0.2cm} \\
	\bottomrule
\end{tabular}\\\vspace{1.5mm}
	\endgroup
	\end{adjustwidth}
Average objective values obtained over 5 runs (due to time constraints, we only used 2 runs for $(P^m)$ and $m=256$) of computing the 2-Wasserstein distance between two uncorrelated Gaussian distributions in $\mathbb{R}^{d}$. For $(P^{m})_{base}$, the parameters are updated taking one infimum update for each supremum update (and we do not include any regularization nor use other techniques for stabilization, such as unrolling or mixtures of generators). For $(P^{m})$, the parameters are updated using 5 unrolling steps of the discriminator (with single updating step for both infimum and supremum) and a mixture of 5 generators. For $(P_\psi)$, we introduce the regularization function $\psi_j^*(x) = \frac{x^2}{150}, j=1,2,$ and take 10 infimum updates for each supremum update. In this case, a single generator is used and no unrolling procedure. The standard deviation of the objective values is computed over the last 5000 iterations.
	\noindent
\end{table}

\begin{figure}[t] 
	\begin{adjustwidth}{-2.5cm}{-2.5cm}
		\centering 
		\includegraphics[width=20cm, height=15cm]{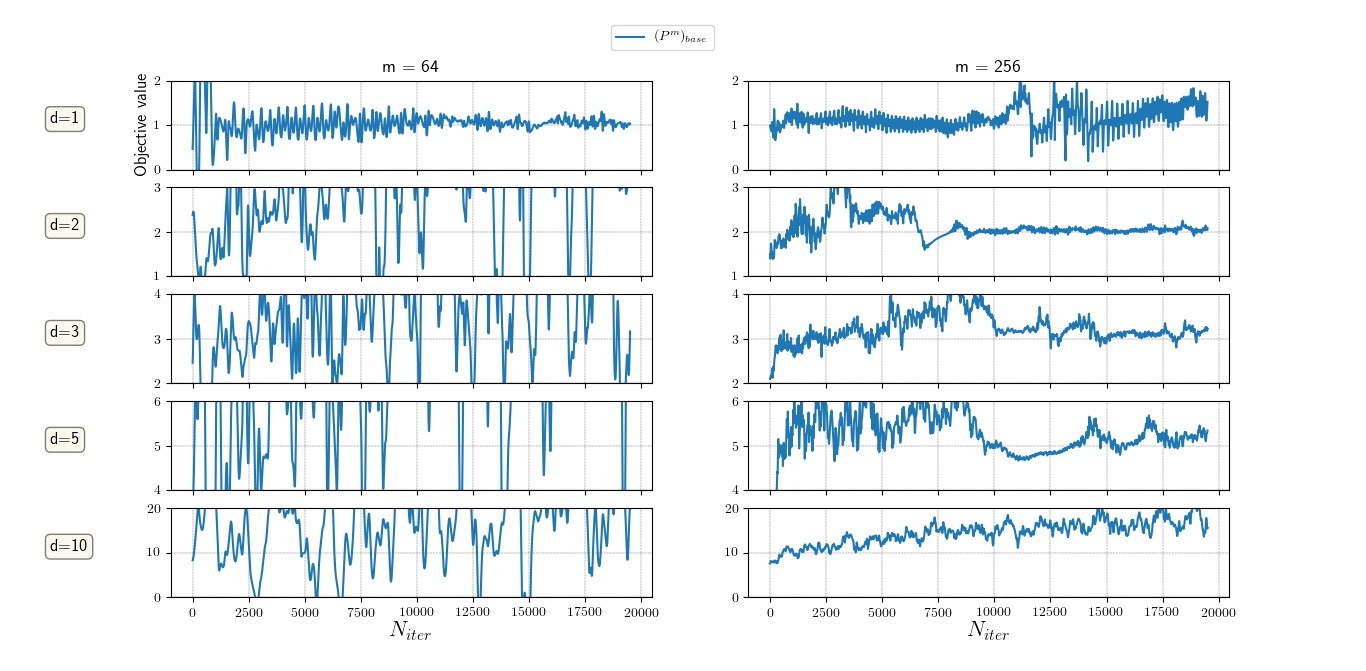} 
	\end{adjustwidth}
	\caption{Numerical convergence observed for the computation of the 2-Wasserstein distance in $\mathbb{R}^{d}$ with formulation $(P^{m})$ and base  optimization procedure (that is, without unrolling and mixture of generators), which we refer to as $(P^{m})_{base}$. The left (resp. right) column shows the convergence when the hidden dimension is set as 64 (resp. 256). The displayed graphs are median values across 5 runs with respect to the standard deviation of the objective values over  the last 5000 iterations.}  \label{fig:W2_2}
\end{figure}

\begin{figure}[t] 
	\begin{adjustwidth}{-2.5cm}{-2.5cm}
	\centering 
	\includegraphics[width=20cm, height=15cm]{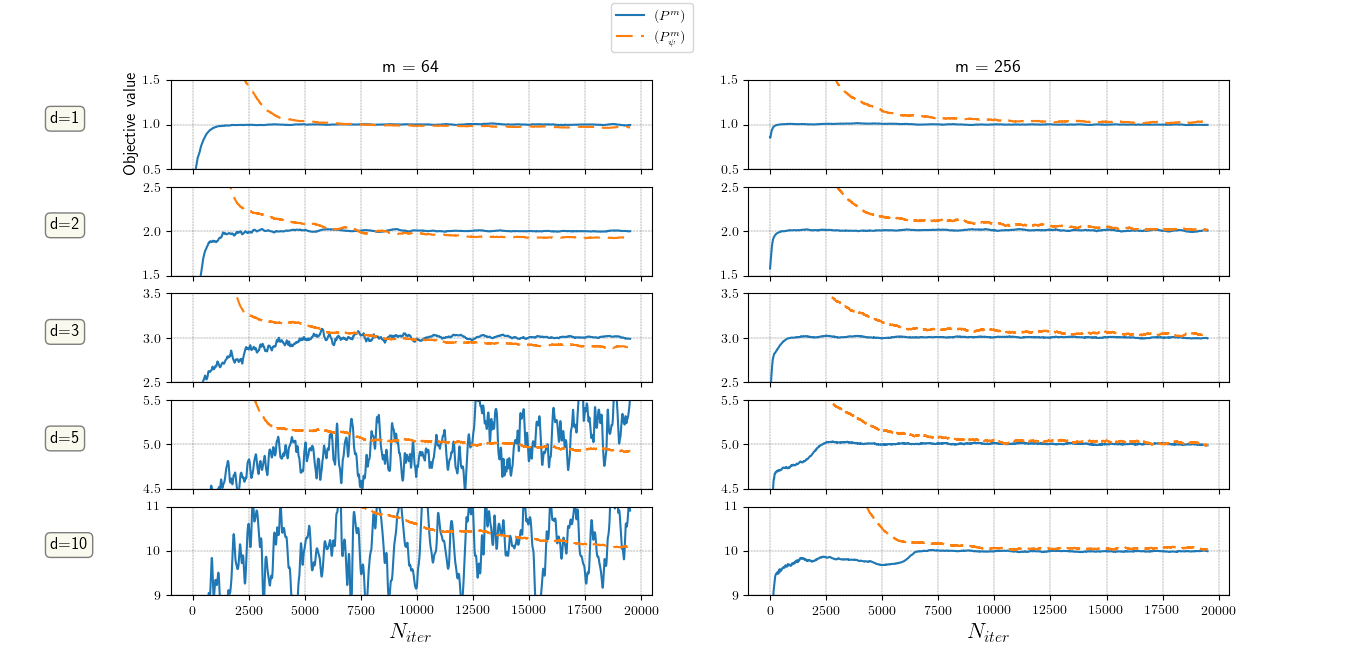} 
	\end{adjustwidth}
	\caption{Comparison of the numerical convergence observed for the computation of the 2-Wasserstein distance in $\mathbb{R}^{d}$ with formulations $(P^{m})$ and $(P^{m}_{\psi})$. The left (resp. right) column shows the convergence when the hidden dimension is set as 64 (resp. 256).   The displayed graphs are median values across 5 runs (2 for $(P^{m})$ and $m = 256$) with respect to the standard deviation of the objective values over  the last 5000 iterations.}  \label{fig:W2}
\end{figure}

The results are provided in Table \ref{W2} and Figures \ref{fig:W2_2} and \ref{fig:W2}. This example corroborates the discussion provided in Section \ref{sec:numerics}. We report three different settings (base case $(P^m_{base})$, combined case $(P^m)$, and $\psi$-regularization $(P^m_\psi)$) for two different network sizes ($m=64$ and $m=256$).
The case $(P^m)_{base}$ results from the simple procedure of using alternating Adam steps for infimum and supremum network, without using regularization, mixtures, or unrolling. 
The case $(P^m)$ corresponds to the combined case from Section \ref{ex:MOT}, i.e., we use both a mixture of 5 generators and 5 steps of unrolling. Finally, the case $(P^m_\psi)$ is the divergence regularization, similar to the one used in Section \ref{ex:distr_constraint}, where we set $\psi_j^*(x) = \frac{x^2}{150}, j=1,2$.

When low computational power is available $(m=64)$, introducing a regularization (formulation $(P^{m}_{\psi})$) helps achieve more stability (even compared to $(P^{m})$), particularly in high-dimensional settings. If, on the other hand, one can increase the hidden dimension ($m=256$) and consequently the runtime, this can also guarantee accuracy and stability of the algorithm both for $(P^{m})$ and $(P^{m}_{\psi})$. The accuracy of $(P^m_{base})$ is limited in either case.

\end{document}